# Fuzzy Limits of Functions

**Mark Burgin**

Department of Mathematics
University of California, Los Angeles
405 Hilgard Ave.
Los Angeles, CA 90095

**Abstract**

The goal of this work is to introduce and study fuzzy limits of functions. Two approaches to fuzzy limits of a function are considered. One is based on the concept of a fuzzy limit of a sequence, while another generalizes the conventional ε-δ definition. It is demonstrated that these constructions are equivalent. Different properties of fuzzy limits of functions are found. Properties of conventional limits are easily deduced from properties of fuzzy limits. In the second part of this work, the concept of fuzzy limits of a function is extended to provide means to define nontrivial continuity of functions on discrete sets. In addition, fuzzy limits of functions are introduced and studied.

## 1. Introduction

Mathematics is an efficient tool for modeling real world phenomena. However, in its essence mathematics is opposite to real world because mathematics is exact, rigorous and abstract while real things and systems are imprecise, vague, and concrete. To lessen this gap, mathematicians elaborate methods that make possible to work with natural vagueness and incompleteness of information using exact mathematical structures. One of the most popular approaches to this problem is fuzzy set theory.

Neoclassical analysis is a synthesis of classical analysis, fuzzy set theory and set valued analysis. In it, ordinary structures of analysis, that is, functions and operators, are studied by means of fuzzy concepts: fuzzy limits, fuzzy continuity, and fuzzy derivatives. For example, continuous functions, which are studied in the classical analysis, become a part of the set of the fuzzy continuous functions studied in neoclassical analysis.

The necessity to launch investigation and implementation of fuzzy principles in the classical analysis, while studying ordinary functions, is caused by several reasons. One of the most important of them is connected with properties of measurements. Any real measurement provides not absolutely precise but only approximate results. For example, it is impossible to find out if any series of numbers obtained in experiments converges or a function determined by measurements is continuous at a given point. Consequently, constructions and methods developed in the classical analysis are only approximations to what exists in reality. In many situations such approximations has been giving a sufficiently adequate representation of studied phenomena. However, scientists and, especially, engineers have discovered many cases in which such methods did not work because classical approach is too rough.

In this paper, we consider basic fuzzy structures that emerge from such fundamental mathematical objects as functions. Analysis of their behavior in computation shows that the classical mathematical idealization is irrelevant for applications. Classical models work well up to some precision. Then they fail. In spite of this, they create an illusion of infinite precision, which does not exist in reality.

Consequently, natural science, as well as computer science, needs new methods to deal with such imprecision. Neoclassical analysis is one of the approaches to deal with imprecision, vagueness, uncertainty and ambiguity of the real world. These new methods and structures provide, in particular, more flexible means for the development of the differential calculus and differential equations for functions defined in discrete spaces than traditional approach based on finite differences.

For simplicity, we consider here only real function on $R$. However, the main constructions and results of this work are valid for a much broader context. It is possible to develop by the same technique and obtain similar results for multidimensional real functions that have the form $f: R^n \to R^m$ and for multidimensional complex functions that have the form $f: C^n \to C^m$.

### Denotations

$N$ is the set of all natural numbers;

$\omega$ is the sequence of all natural numbers;

$Z$ is the set of all integer numbers;

$\varnothing$ is the empty set;

$R$ is the set of all real numbers;

$R^+$ is the set of all non-negative real numbers;

$R^{++}$ is the set of all positive real numbers;

$R_\infty = R \cup \{\infty, -\infty\}$;

if $a$ is a real number, then $|a|$ or $\|a\|$ denotes its absolute value or modulus;

if $a$ is a real number and $t \in R^{++}$, then $O_t a = \{ x \in R; a - t < x < a + t \}$ is a neighborhood of $a$;

$\rho(x,y) = | x - y |$ for $x, y \in R$;

if $U \subseteq R$, then the diameter $D(U)$ of the set $U$ is equal to $\sup\{ \rho(x, y); x, y \in U \}$;

if $l = \{a_i \in M; i \in \omega\}$ is a sequence, and $f: M \to L$ is a mapping, then $f(l) = \{f(a_i); i \in \omega\}$;

the logical symbol $\forall$ means "for any";

the logical symbol $\exists$ means "there exists";

if $l = \{a_i \in M; i \in \omega\}$ is a sequence, and $f: M \to L$ is a mapping, then $f(l) = \{f(a_i); i \in \omega\}$;

$a = r$-lim $l$ means that a number $a$ is an $r$-limit of a sequence $l$;

$a = (B, r)$-lim $l$ means that a number $a$ is a conditional $r$-limit or $(B, r)$-limit of a sequence $l$;

if $U$ is a correspondence of a set $X$ to a set $Y$ (binary relation on $X$ ), then $U(x) = \{y \in Y; (x, y) \in U \}$ and $U^{-1}(y) = \{x \in X; (x, y) \in U \}$;

If $f: \boldsymbol{R} \to \boldsymbol{R}$ is a partial function, then Dom $f$ is the domain and Rg $f$ is the range of $f$.

## 2. Fuzzy limits of sequences

The theory of fuzzy limits of functions is based on the theory of fuzzy limits of sequences (Burgin, 2000). That is why we begin our exposition with the main concepts and constructions from the latter theory.

Let $r \in \boldsymbol{R}^+$ and $l = \{a_i \in \boldsymbol{R}; i \in \omega\}$ be a sequence of real numbers.

**Definition 2.1.** A number $a$ is called an $r$-limit of a sequence $l$ (it is denoted by $a = r\text{-lim}_{i \to \infty} a_i$ or $a = r$-lim $l$ ) if for any $\varepsilon \in \boldsymbol{R}^{++}$ the inequality $\rho(a, a_i) < r + \varepsilon$ is valid for almost all $a_i$, i.e., there is such n that for any $i > n$, we have $\rho(a, a_i) < r + \varepsilon$.

In this case, we say that $l$ $r$-converges to $a$ and denote it by $l \to_q a$.

**Example 2.1.** Let $l = \{1/i; i \in \omega\}$. Then 1 and -1 are 1-limits of $l$ ; 1/2 is a (1/2)-limit of $l$ but 1 is not a (1/2)-limit of $l$.

Informally, $a$ is an $r$-limit of a sequence $l$ for an arbitrarily small $\varepsilon$ the distance between $a$ and all but a finite number of elements from $l$ is smaller than $r + \varepsilon$ .

When $r = 0$, the $r$-limit coincides with the conventional limit of a sequence. As a result, the concept of an $r$-limit is a natural extension of the concept of conventional limit. However, the concept of an $r$-limit actually extends the conventional construction of a limit (cf., Burgin, 2000).

**Lemma 2.1.** If $a = r$-lim $l$, then $a = q$-lim $l$ for any $q > r$.

**Definition 2.2.** a) A number $a$ is called a fuzzy limit of a sequence $l$ if it is an $r$-limit of $l$ for some $r \in \mathbf{R}^+$.

b) a sequence $l$ fuzzy converges if it has a fuzzy limit.

**Example 2.2.** Let us consider sequences $l = \{1 + 1/i;\ i \in \omega\}$, $h = \{1 + (-1)^i;\ i \in \omega\}$, and $k = \{1 + [(1 - i)/i]^i\ ;\ i \in \omega\}$. Sequence $l$ has the conventional limit equal to 1 and many fuzzy limits (e.g., 0, 0.5, 2 are 1-limits of $l$ ). Sequence $h$ does not have the conventional limit but has different fuzzy limits (e.g., 0 is a 1-limit of $h$, while 1, -1, and 1/2 are 2-limits of $h$). Sequence $k$ does not have the conventional limit but has a variety of fuzzy limits (e.g., 1 is a 1-limit of $k$, while 2, 0, 1.5, 1.7, and 0.5 are 2-limits of $k$ ).

Thus, we see that many sequences that do not have the conventional limit have a quantity of fuzzy limits.

**Lemma 2.2.** If $a = q$-lim $h$, then $a = q$-lim $k$ for any subsequence $k$ of $h$.

**Remark 2.1.** The measure of convergence of $l$ to points from $\mathbf{R}$ defines the normal fuzzy set Lim $l = [L, \mu(x = \lim l)]$ of fuzzy limits of $l$ (Burgin, 2000).

Let $l = \{a_i \in \mathbf{R};\ i \in \omega\}$ be a bounded sequence.

**Lemma 2.3.** If for any convergent subsequence $k$ of $l$, we have $a = r$-lim $k$, then $a = r$-lim $l$.

Proof. Let us assume that the condition of the lemma is satisfied, but $a \neq r$-lim $l$. Then there is $\varepsilon \in \mathbf{R}^{++}$ such that for infinitely many elements $a_i$, we have $\rho(a, a_i) > r + \varepsilon$. Let us take these elements $a_{i1}$, $a_{i2}$, ... , $a_{in}$, ... . According to our assumption, the sequence $h = \{a_{in} \in \mathbf{R};\ n \in \omega\}$ is bounded, as a subsequence of a bounded sequence. Consequently, $h$ has a convergent subsequence $h = \{a_{in} \in \mathbf{R};\ n \in \omega\}$ (cf., for example, (Randolph, 1968)). If $d = \text{lim } k$, then $\rho(a, d) \geq r + \varepsilon$ (cf., for example, (Burgin, 2000)). According to Definition 2.1, the point $a$ is not an $r$-limit of the sequence $l$. This contradicts our assumption and by the principle of excluded middle, concludes the proof.

**Definition 2.3.** $\infty$ ($-\infty$) is an $r$-limit of $l$ if almost all elements $a_i$ are bigger (less) than $r$ ($-r$).

**Definition 2.4.** A number $b$ is called a weak or partial limit (weak $r$-limit) of a sequence $l = \{a_i \in R; i \in \omega\}$ (it is denoted by $b =$ wlim $l$ or $b = r$-wlim $l$) if there is a subsequence $h$ the sequence $l$ such that $a = \lim h$ ($b = r$-wlim $l$).

Weak limits of sequences are used in the theory of automata working with infinite words. These automata are called timed and reflect temporal behavior of finite-state reactive programs such as communication and synchronization protocols. Examples of timed automata are Büchi (1960) and Muller (1963) automata.

**Definition 2.5.** A Büchi automaton $A$ is a finite automaton input of which is an infinite sequence $l = \{a_i \in \Sigma; i \in \omega\}$ of symbols. Each input sequence $l$ generates a sequence $\sigma(l)$ of states of $A$ that is called a run of $A$. As it is usual for finite automata, a subset $F$ of the set $Q_A$ of states of $A$ is chosen, elements of which are called finite or accepting states. An input sequence $l$ is accepted if some states from $F$ occur infinitely often in $\sigma(l)$.

**Definition 2.6.** A Muller automaton $A$ is a finite automaton input of which is an infinite sequence $l = \{a_i \in \Sigma; i \in \omega\}$ of symbols. Each input sequence $l$ generates a sequence $\sigma(l)$ of states of $A$ that is called a run of $A$. In contrast to Büchi automata, a set $F$ of subsets of the set $Q_A$ of states of $A$ is chosen. Elements of $F$ are called finite or accepting sets. An input sequence $l$ is accepted by $A$ if the set of all states repeating infinitely often in $\sigma(l)$ belongs to $F$.

Informally, the accepting states of these automata designate states through which a run must pass infinitely often. More formally, on an infinite word, every run of an infinite automaton will visit some set of states infinitely often. An infinite automaton accepts a word iff this set intersects the set F of accepting states. In such a way, Büchi and Muller automata define ω–regular languages. While accepting a word, a run need not visit every accepting state infinitely often; it only needs visit some accepting state infinitely often. It means that accepting states have to be weak limits of all states in an automaton run. It is made precise and formalized by the following proposition.

Let the set $Q_A$ is considered as a topological space with the discrete topology.

**Proposition 2.1.** a) An input sequence $l$ is accepted by a Büchi automaton $A$ with a set $F$ of final states if and only if $\sigma(l)$ has a weak limit in $F$.

b) An input sequence *l* is accepted by a Muller automaton *A* with a set *F* of accepting sets if and only if the set of all weak limits of σ(*l*) belongs to *F*.

Timed automata provide the foundation for theories and techniques for system analysis techniques such as computer-aided verification. Using timed automata, we can model desired properties of systems and protocols, such as "whenever the system receives a request, it eventually produces an acknowledgement." Then, given a model M of a system, we can see whether every run of the system satisfies the property by checking whether the language produced by M is accepted by the automaton for the property.

The untiming construction for timed automata forms the basis for verification algorithms for the branching temporal logics, as well as for realizability checking and synthesis of reactive software modules (Vardi and Wolper, 1986; 1994).

Neoclassical analysis makes possible not only to extend ordinary concepts obtaining new results for classical structures, but also provides for elaboration of new useful concepts. One of such concepts is given in the definition of fuzzy limits of sets of sequences.

**Definition 2.7.** A number *a* is called an *r*-limit of a set E = { $l_j$ ; j∈ω } of sequences of real numbers (it is denoted by *a* = *r*-lim E ) if *a* is an *r*-limit of each sequence $l_i$ from E.

**Remark 2.2.** If E has a 0-limit *a*, then this 0-limit is unique and all sequences from E converge to *a*. In contrast to this sequences from a given set may have different limits but a common fuzzy limit. It is demonstrated by the following example.

**Example 2.3.** Let us consider the set E = { {$1/2^n$; *n* = 1, 2,...}, {1+ $1/3^n$; *n* = 1, 2,...}, {2+ $1/5^n$; *n* = 1, 2,...} }. It has a 1-limit 1, but it does not have a limit because the first sequence converges to 0, the second sequence converges to 1, and the third sequence converges to 2.

**Theorem 2.1 (Reduction Theorem).** If the set E of sequences is finite or countable, then there is a sequence *l* = { $a_i$ ∈ **R**; *i*∈ω} such that for any *a* ∈ **R** and any *r* ∈ **R**$^+$, we have *a* = *r*-lim E if and only if *a* = *r*-lim *l*.

## 3. Fuzzy limits of functions

At first, we develop and study a construction similar to fuzzy limits of sequences because it is only one level of generality higher than the classical concept of a limit of a function. This construction allows one to model and study continuous mappings from continuous to discrete spaces.

Let $r \in \mathbf{R}^+$ and $f: \mathbf{R} \to \mathbf{R}$ be a partial function.

**Definition 3.1**. A number $b$ is called an $r$-limit of a function $f$ at a point $a \in \mathbf{R}_\infty$ (it is denoted by $b = r\text{-}\lim_{x \to a} f(x)$) if for any sequence $l = \{a_i \in \text{Dom}\, f;\ i \in \omega,\ a_i \neq a\}$, the condition $a = \lim l$ implies $b = r\text{-}\lim_{i \to \infty} f(a_i)$.

**Remark 3.1.** It is possible to define an $r$-limit of a function $f$ at a point $a \in \mathbf{R}$ independently of the concept of an $r$-limit of a sequence (cf. Theorem 3.1). This allows one to get the latter concept as a particular case of the former concept. Really, a sequence with values in $\mathbf{R}$ is a partial function $f: \mathbf{R} \to \mathbf{R}$ that is defined only for natural numbers, i.e., $f: \mathbf{N} \to \mathbf{R}$. Then an $r$-limit of a sequence $l = \{a_i \in \mathbf{R};\ i \in \omega\}$ is an $r$-limit of the function $f$ at the point $\infty$.

**Remark 3.2.** The concept of an $r$-limit of a function allows one to develop a theory of limit for functions that take values in discrete sets.

**Lemma 3.1.** If $\text{Dom}\, f = \mathbf{R}$, then $b = 0\text{-}\lim_{x \to a} f(x)$ if and only if $b = \lim_{x \to a} f(x)$ in the classical sense.

This result demonstrates that the concept of an $r$-limit of a function is a natural extension of the concept of conventional limit of a function. However, the concept of an $r$-limit actually extends the conventional construction of a limit (cf. Example 3.1).

**Remark 3.3.** The condition $\text{Dom}\, f = \mathbf{R}$ is essential in Lemma 3.1 because if $\text{Dom}\, f$ is a discrete set, then all sequences converging to any element from $\text{Dom}\, f$ stabilize after some number of its elements. This results in the property that any function $f$ defined on a discrete set has the limit at any point of this set, which is equal to value of $f$, and thus, $f$ is continuous.

**Lemma 3.2.** If $b = r\text{-}\lim_{x \to a} f(x)$, then $b = q\text{-}\lim_{x \to a} f(x)$ for any $q > r$.

**Theorem 3.1.** The condition $b = r\text{-}\lim_{x \to a} f(x)$ is valid if and only if for any open neighborhood $Ob$ of $b$ that contains the interval $[b - r, b + r]$ there is a neighborhood $Oa$ of $a$ such that $f(Oa \cap \text{Dom } f) \subseteq Ob$.

Proof. Necessity. Let $b = r\text{-}\lim_{x \to a} f(x)$ and $Ob$ is an open neighborhood of $b$ that contains the interval $[b - r, b + r]$. Let us suppose that in any neighborhood $Oa$ of $a$, there is a point $a_i \neq a$ such that $f(a_i)$ does not belong to $Ob$. We can take a sequence of neighborhoods $O_i a$ such that $O_i a \subseteq O_{i-1} a$ for all $i = 2, 3, \ldots$ , in any neighborhood $O_i a$, there is a point $a_i \neq a$ such that $f(a_i)$ does not belong to $Ob$, and $\cap O_i a = \{a\}$. Then we have a sequence

$$l = \{a_i \in R; i \in \omega, a_i \neq a\},$$

for which the condition $a = \lim l$ does not imply $b = r\text{-}\lim_{x \to a} f(a_i)$. This contradicts the initial condition that $b = r\text{-}\lim_{x \to a} f(x)$ and concludes the proof of necessity.

Sufficiency. Let for any open neighborhood $Ob$ of $b$ that contains the interval $[b - r, b + r]$ there is a neighborhood $Oa$ of $a$ such that $f(Oa) \subseteq Ob$ and $l = \{a_i \in R; i \in \omega, a_i \neq a\}$ is a sequence such that $a = \lim l$. Then almost all elements of $l$ belong to $Oa$. Consequently, almost all elements $f(a_i)$ belong to $Ob$. By the definition of $r$-limit, we have $b = r\text{-}\lim_{i \to \infty} f(a_i)$. As the sequence $l$ is chosen arbitrarily, by Definition 3.1, $b = r\text{-}\lim_{x \to a} f(x)$. This concludes the proof of the theorem.

Theorem 2.1 from (Burgin, 2000) implies the following result.

**Theorem 3.2.** If $b = r\text{-}\lim_{x \to a} f(x)$ and $b > d + r$, then there is a neighborhood $Oa$ of $a$ such that $f(x) > d$ for all $x$ from $Oa \cap \text{Dom } f$.

**Corollary 3.1.** If $b = r\text{-}\lim_{x \to a} f(x)$ and $b > d + r$, then for any sequence $l = \{a_i \in \text{Dom } f; i \in \omega, a_i \neq a\}$ with $a = \lim l$, we have $f(a_i) > d$ for almost all $a_i$ from $l$.

**Corollary 3.2.** If $f(x) \leq q$ for all $x$ from some neighborhood $Oa$ of $a$ and $a = r\text{-}\lim_{x \to a} f(x)$, then $a \leq q + r$.

These results allow us to obtain important properties of conventional limits of functions that are studied in courses of calculus.

**Corollary 3.3**. If $a = \lim_{x \to a} f(x)$ and $a > b$, then $f(x) > b$ for all $x$ from some neighborhood $Oa$ of $a$.

**Corollary 3.4**. If $a = \lim_{x \to a} f(x)$ and $a > 0$, then $f(x) > 0$ for all $x$ from some neighborhood $Oa$ of $a$.

**Corollary 3.5**. If $f(x) \leq q$ for all $x$ from some neighborhood $Oa$ of $a$ and $a = \lim_{x \to a} f(x)$, then $a \leq q$.

**Definition 3.2**. a) A number $a$ is called a *fuzzy limit* of a function $f(x)$ at a point $a \in \mathbf{R}$ if it is an $r$-limit of a function $f(x)$ at a point $a$ for some $r \in \mathbf{R}^+$.

b) a function $f(x)$ *fuzzy converges* at a point $a \in \mathbf{R}$ if it has a fuzzy limit at this point.

**Example 3.1.** Let us consider the following function $f(x)$:

$$f_n(x) = \begin{cases} 1 + 1/i & \text{when } x = 1 - 1/i \text{ for } i = 1, 2, 3, \ldots; \\ 2 - 1/i & \text{when } x = 1 + 1/2i \text{ for } i = 1, 2, 3, \ldots; \\ 1 + (-1)^i & \text{when } x = 1 + 1/(2i+1) \text{ for } i = 1, 2, 3, \ldots; \\ x & \text{otherwise}. \end{cases}$$

The function $f(x)$ does not have the conventional limit at the point 1 but has different fuzzy limits at this point. For instance, 1 is a 1-limit of $f$, while 2, 0, 1.5, 1.7, and 0.5 are 2-limits of $f$.

Thus, we see that many functions that do not have the conventional limit at some point have lots of fuzzy limits.

**Remark 3.4.** The measure of convergence of a function $f(x)$ at a point $a$ to points from $\mathbf{R}$ defines the normal fuzzy set $\text{Lim}_a f = [L, \mu(x = \lim_{x \to a} f(x))]$ of fuzzy limits of the function $f$ at the point $a$.

**Remark 3.5.** The property of fuzzy limits of a function $f$ at a point $a$ implies that for a bounded function $f$ at the point $a$ all real numbers become its fuzzy limits at this points. However, the advantage of our approach is that we classify all these fuzzy limits by the measure of convergence. Thus, it is not reasonable to consider the set of all fuzzy limits of a bounded function *f* at a point *a*, while the fuzzy set Lim$_a$ *f* of all fuzzy limits of a bounded function *f* at a point *a* gives a lot of information about this function.

**Definition 3.3.** A number *b* is called a *weak r-limit* of a function *f(x)* at a point $a \in \mathbf{R}$ (it is denoted by $b = r\text{-wlim}_{x \to a} f(x)$) if there is a sequence $l = \{a_i \in \mathbf{R}; i \in \omega, a_i \neq a\}$ such that $a = \lim l$ and $b = r\text{-}\lim_{i \to \infty} f(a_i)$.

The concept of a weak *r*-limit of a function gives a new concept for the classical case.

**Definition 3.4.** A number *b* is called a *weak limit* of a function *f(x)* at a point $a \in \mathbf{R}$ (it is denoted by $b = \text{wlim}_{x \to a} f(x)$) if there is a sequence $l = \{a_i \in \text{Dom} f ; i \in \omega, a_i \neq a\}$ such that $a = \lim l$ and $b = \lim_{i \to \infty} f(a_i)$.

In the classical mathematical analysis, some cases of weak limits (such as $\underline{\lim} = \lim \inf$ or $\overline{\lim} = \lim \sup$ ) are considered. A general case of partial limits is treated in (Randolph, 1968) where weak limits, i.e., weak 0-limits by Lemma 3.1, are called subsequential limits.

**Remark 3.6.** In the theory of hypernumbers and extrafunctions, weak limits are used to build the spectrum of a hypernumber (Burgin, 2002).

**Lemma 3.3.** If Dom $f = \mathbf{R}$, then the number $b = \text{wlim}_{x \to a} f(x)$ if and only if $b = 0\text{-wlim}_{x \to a} f(x)$.

**Lemma 3.4.** If $b = r\text{-wlim}_{x \to a} f(x)$, then $b = q\text{-wlim}_{x \to a} f(x)$ for any $q > r$.

**Definition 3.5**. A number *a* is called a *weak fuzzy limit* of a function *f* at a point $a \in \mathbf{R}$ if it is a weak *r*-limit of a function *f* at the point *a* for some $r \in \mathbf{R}^+$.

Fuzzy limits and weak fuzzy limits of functions are invariant with respect to a continuous monotonous change of the scale of the argument.

**Theorem 3.3.** If $b = r\text{-}\lim_{x \to a} f(x)$ ($b = r\text{-wlim}_{x \to a} f(x)$) ) and *g(x)* is a monotonous continuous function with the continuous inverse function, then $b = r\text{-}\lim_{x \to h(a)} f(g(x))$ ( $b = r\text{-wlim}_{x \to h(a)} f(g(x))$ ) where $h(a) = g^{-1}(x)$.

<u>Proof</u>. (a) Let us assume that $b = r\text{-}\lim_{x \to a} f(x)$ and take some sequence $l = \{a_i \in \mathbf{R}; i \in \omega\}$ that satisfies the condition $h(a) = \lim_{i \to \infty} a_i$. As *g(x)* is a continuous function, we have $a = g(h(a)) = g(g^{-1}(a)) = \lim_{i\to\infty} g(a_i)$. Then by Definition 3.1, we have $b = r\text{-}\lim_{x\to h(a)} f(g(a_i))$. Thus, $b = r\text{-}\lim_{x\to h(a)} f(g(x))$ because $l$ is an arbitrary sequence.

(b) Let us assume that $b = r\text{-}\lim_{x\to a} f(x)$ and take such sequence $l = \{a_i \in R; i \in \omega\}$ that satisfies the condition $h(a) = \lim_{i\to\infty} a_i$. By Definition 3.2, such a sequence exists. Then $h(a) = \lim_{i\to\infty} g^{-1}(a_i)$ because $g^{-1}(x)$ is a continuous function. By the definition of an $r$-limit, this implies that for any $\varepsilon > 0$ there is such a $\delta > 0$ that whenever $|g^{-1}(a_i) - h(a)| \le \delta$, we have $|f(a_i) - a| = |f(g(g^{-1}(a_i))) - a| \le \varepsilon$. Consequently, $b = r\text{-w}\lim_{x\to h(a)} f(g(x))$.

Theorem is proved.

**Corollary 3.6.** If $b = r\text{-}\lim_{x\to a} f(x)$ ($b = r\text{-w}\lim_{x\to a} f(x)$) and $k \in R^+$, then $b = r\text{-}\lim_{x\to ha} f(kx)$ ($b = r\text{-w}\lim_{x\to ha} f(kx)$) where $h = k^{-1}$.

Theorem 2.2 from (Burgin, 2000) implies the following result.

**Theorem 3.4**. For an arbitrary number $r \in R^+$, all $r$-limits at a point $a$ of a locally bounded function $f$ belong to some finite interval, the length of which is equal to $2r$.

**Remark 3.7.** For unbounded functions, this property can be false.

Theorem 3.4 and Lemma 3.1 imply the following result.

**Corollary 3.7** (Any course of calculus, cf., for example, (Ross, 1996)). If the limit of a function at some point exists, then it is unique.

**Theorem 3.5.** Let $b = r\text{-}\lim_{x\to a} f(x)$ and $c = q\text{-}\lim_{x\to a} g(x)$. Then:

a) $b + c = (r+q)\text{-}\lim_{x\to a} (f + g)(x)$;

b) $b - c = (r+q)\text{-}\lim_{x\to a} (f - g)(x)$;

c) $kb = (|k|\cdot r)\text{-}\lim_{x\to a} (kf)(x)$ for any $k \in R$ where $(k f)(x) = k \cdot f(x)$.

This theorem is derived as a corollary from a more general Theorem 3.11 proved below.

**Corollary 3.8** (Any course of calculus, cf., for example, (Ross, 1996)). Let $b = \lim_{x\to a} f(x)$ and $c = \lim_{x\to a} g(x)$. Then:

a) $b + a = \lim_{x\to a} (f + g)(x)$;

b) $b - c = \lim_{x\to a} (f - g)(x)$;

c) $ka = \lim_{x\to a} (kf)(x)$ for any $k \in R$.

**Definition 3.6.** A function $f: \mathbf{R} \to \mathbf{R}$ is *bounded* at a point $a \in \mathbf{R}$ if there exists a number $q$ and a neighborhood $Oa$ of the point $a$ such that for any $x$ from $Oa$ the inequality $\rho(f(a), f(x)) < q$ is valid.

**Theorem 3.6.** A function $f(x)$ fuzzy converges at a point $a$ if and only if $f(x)$ is bounded at this point.

A proof of this statement is based on connections between fuzzy continuous functions and fuzzy limit and utilizes a corresponding result for fuzzy continuous functions (cf., (Burgin, 1995)).

As in the case of sequences, fuzzy limits of functions are more adequate to real situations and provide means for more relevant models than classical limits. However, to treat functions defined on discrete sets, we need a more advanced concept of a limit.

Let $q, r \in \mathbf{R}^+$ and $f: \mathbf{R} \to \mathbf{R}$ be a partial function.

**Definition 3.7.** A number $b$ is called a $(q, r)$-*limit* of a function $f$ at a point $a \in \mathbf{R}_\infty$ (it is denoted by $b = (q, r)\text{-}\lim_{x \to_q a} f(x)$) if for any sequence $l = \{a_i \in \text{Dom } f; i \in \omega, a_i \neq a\}$, the condition $a = q\text{-}\lim l$ implies $b = r\text{-}\lim_{i \to \infty} f(a_i)$.

**Example 3.2.** Let us consider the function $f(x) = x/|x|$. The graph is presented in Figure 1. The classical limit $\lim_{x \to 0} (x/|x|)$ does not exist. At the same time, this function has fuzzy and fuzzy fuzzy limits. For instance, $0 = 1\text{-}\lim_{x \to 0} (x/|x|)$, $0 = (0, 1)\text{-}\lim_{x \to 0} (x/|x|)$ and $1 = (1, 2)\text{-}\lim_{x \to_2 0} (x/|x|)$.

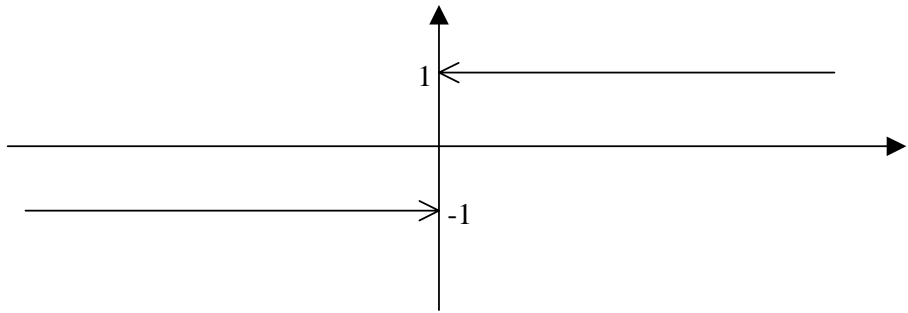

**Fig 1.**

**Remark 3.8.** It is possible to define a $(q, r)$-limit of a function $f$ at a point $a \in \mathbf{R}$ independently of the concept of an $r$-limit of a sequence (cf. Theorem 3.7).

**Lemma 3.5.** The number $b = r\text{-}\lim_{x \to a} f(x)$ if and only if $b = (0, r)\text{-}\lim_{x \to a} f(x)$.

This result demonstrates that the concept of a $(q, r)$-limit of a function is a natural extension of the concept of an $r$-limit of a function.

Lemma 3.2 and Definition 3.7 imply the following result.

**Lemma 3.6.** If $b = (q, r)\text{-}\lim_{x \to^q a} f(x)$, then $b = (u, v)\text{-}\lim_{x \to^u a} f(x)$ for any $v \geq r$ and $u \leq q$.

Indeed, let us take some sequence $l = \{a_i \in \text{Dom} f; i \in \omega, a_i \neq a\}$ such that $a = q\text{-}\lim l$. As $u \leq q$, for any sequence $l$, the condition $a = u\text{-}\lim l$ implies (by Lemma 2.1) the condition $a = q\text{-}\lim l$. If $b = (q, r)\text{-}\lim_{x \to a} f(x)$, then by Definition 3.7, for any sequence $h = \{c_i \in \text{Dom} f; i \in \omega, c_i \neq a\}$, the condition $a = q\text{-}\lim l$ implies $b = r\text{-}\lim_{i \to \infty} f(c_i)$. In particular, we have $b = r\text{-}\lim_{i \to \infty} f(a_i)$. As $v \geq r$, the condition $b = r\text{-}\lim_{i \to \infty} f(a_i)$ implies (by Lemma 2.1) the condition $b = v\text{-}\lim_{i \to \infty} f(a_i)$. As $l$ is an arbitrary sequence of elements from the domain $\text{Dom} f$, $b = (u, v)\text{-}\lim_{x \to a} f(x)$ by Definition 3.7.

**Theorem 3.7.** The condition $b = (q, r)\text{-}\lim_{x \to a} f(x)$ is valid if and only if for any $x \in \text{Dom} f$ and any $\varepsilon \in \mathbf{R}^{++}$ there is $\delta \in \mathbf{R}^{++}$ such that the inequality $\rho(x, a) < r + \delta$ implies the inequality $\rho(f(x), b) < r + \varepsilon$.

Proof. Necessity. Let $b = (q, r)\text{-}\lim_{x \to^q a} f(x)$ and $Ob$ is an open neighborhood of $b$ that contains the interval $[b - r, b + r]$. Let us suppose that in any neighborhood $Oa$ of $a$ such that $Oa$ contains the interval $[a - q, a + q]$, there is a point $a_i \neq a$ such that $f(a_i)$ does not belong to $Ob$. We can take a sequence of neighborhoods $O_i a$ such that $O_i a \subseteq O_{i-1} a$ for all $i = 2, 3, \ldots$, in any neighborhood $O_i a$, there is a point $a_i \neq a$ such that $f(a_i)$ does not belong to $Ob$, and $\cap O_i a = \{a\}$. Then we have a sequence

$$l = \{a_i \in \mathbf{R}; i \in \omega, a_i \neq a\},$$

for which the condition $a = q\text{-}\lim l$ does not imply $b = (q, r)\text{-}\lim_{i \to \infty} f(a_i)$. This contradicts the initial condition that $b = (q, r)\text{-}\lim_{x \to a} f(x)$ and concludes the proof of necessity.

Sufficiency. Let for any open neighborhood $Ob$ of $b$ that contains the interval $[b - r, b + r]$ there is a neighborhood $Oa$ of $a$ such that $Oa$ contains the interval $[a - q, a + q]$, $f(Oa) \subseteq Ob$, and $l = \{a_i \in R; i \in \omega, a_i \neq a\}$ is a sequence such that $a = q\text{-lim } l$. Then almost all elements of $l$ belong to $Oa$. Consequently, almost all elements $f(a_i)$ belong to $Ob$. By the definition of $r$-limit, we have $b = r\text{-lim}_{i \to \infty} f(a_i)$. As the sequence $l$ is chosen arbitrarily, by Definition 3.6, $b = (q, r)\text{-lim}_{x \to a} f(x)$. This concludes the proof of the theorem.

Theorem 3.7 and Lemma 3.6 imply the following result.

**Theorem 3.8.** If $b = (q, r)\text{-lim}_{x \to^q a} f(x)$ and $b > d + r$, then there is a neighborhood $Oa$ of $a$ such that $Oa$ contains the interval $[a - q, a + q]$ and $f(x) > d$ for all $x$ from $Oa \cap \text{Dom } f$.

**Corollary 3.9.** If $b = (q, r)\text{-lim}_{x \to^q a} f(x)$ and $b > d + r$, then for any sequence $l = \{a_i \in R; i \in \omega, a_i \neq a\}$ with $a = r\text{-lim } l$, we have $f(a_i) > d$ for almost all $a_i$ from $l$.

**Corollary 3.10.** If $f(x) \leq d$ for all $x$ from some neighborhood $Oa$ of $a$ that contains the interval $[a - q, a + q]$ and $b = (q, r)\text{-lim}_{x \to^q a} f(x)$, then $b \leq d + r$.

**Definition 3.8**. a) A number $a$ is called a *fuzzy-fuzzy limit* of a function $f(x)$ at a point $a \in R$ if it is an $(q, r)$-limit of a function $f$ at the point $a$ for some $q, r \in R^+$.

b) a function $f(x)$ *fuzzy-fuzzy converges* at a point $a \in R$ if it has a fuzzy-fuzzy limit at this point.

**Example 3.3.** Let us consider the function $g(x) = [x]$ where $[x]$ is the largest integer that is less than or equal to $x$. The graph of $g(x)$ is presented in Figure 2. The classical limit $\lim_{x \to n} [x]$ does not exist for any integer number $n$ because $\lim_{x \to n-} [x] = n - 1$ and $\lim_{x \to n+} [x] = n$. At the same time, this function has fuzzy and fuzzy-fuzzy limits. For instance, $2.5 = (0, \frac{1}{2})\text{-lim}_{x \to 3} [x]$, $5 = (0.9, 1)\text{-lim}_{x \to^{0.9} 5} [x]$, and $3 = (1, 2)\text{-lim}_{x \to^1 3} [x]$. Thus, 2.5 and 3 are fuzzy-fuzzy limits of $g(x)$.

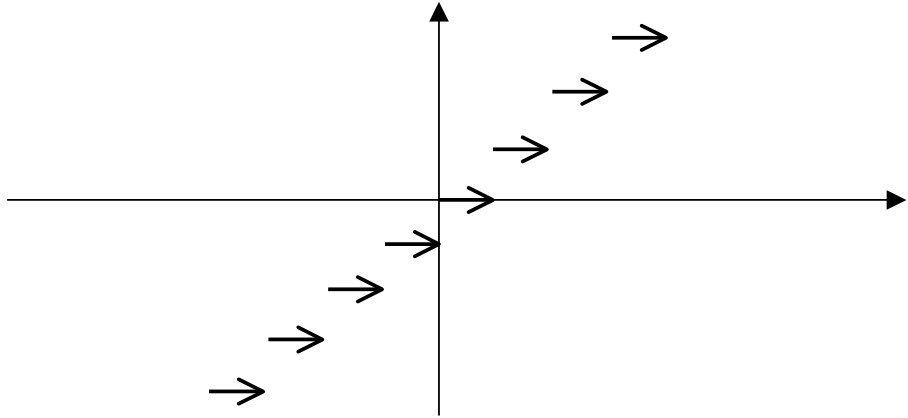

**Fig 2.**

**Remark 3.9.** The concept of a fuzzy limit takes into account fuzziness only in one dimension. Taking the conventional *x-y* coordinates, this is the *y*- dimension. The concept of a fuzzy-fuzzy limit takes into account fuzziness in two dimensions. Namely, they are both *x*- and *y*-dimensions.

**Remark 3.10.** The concept of an ($q$, $r$)-limit of a function allows one to develop a theory of limit for functions that are defined and take values in discrete sets. Let us consider some examples.

**Example 3.4.** Let us consider the function $f(x) = x$ defined for all integers. The classical limit $\lim_{x \to n} x$, as well as any *r*-limit *r*-$\lim_{x \to n} x$, becomes meaningless (being automatically equal to *n* for any integer number *n*) for discerning continuous, in some sense, functions on integers. The cause is that *x* can never come infinitely close to *n* without being equal to *n*. At the same time, this function has different fuzzy-fuzzy limits that provide for a broader concept of continuity on discrete sets. For instance, $2.5 = (1, 1.5)$-$\lim_{x \to^1 3} [x]$, $5 = (0.5, 1)$-$\lim_{x \to^{0.5} 5} [x]$, and $3 = (1, 1)$-$\lim_{x \to^1 3} [x]$. Thus, 2.5 and 3 are fuzzy-fuzzy limits of $f(x)$.

**Remark 3.11.** Fuzzy limits or *r*-limits of functions form a base for the concept of fuzzy continuity or *r*-continuity (cf. Burgin, 1995; 2004). In a similar way, fuzzy-fuzzy limits or $(q, r)$-limits of functions form a base for the concept of fuzzy-fuzzy continuity or $(q, r)$-continuity.

The concept of an $(q, r)$-limit of a function allows one to develop a theory of limit for functions that are defined and take values in discrete sets.

There is another way to construct fuzzy-fuzzy limits of functions.

Let $r \in \mathbf{R}^+$, $f: \mathbf{R} \to \mathbf{R}$ be a partial function, and $[u, v]$ be an interval in $\mathbf{R}$.

**Definition 3.9**. A number $b$ is called an $r$-limit of a function $f$ at an interval $[u, v]$ (it is denoted by $b = r\text{-lim}_{x \to [u, v]} f(x)$) if for any point $a$ from the interval $[u, v]$ and any sequence $l = \{a_i \in \text{Dom} f; i \in \omega, a_i \neq a\}$ with $a = \lim l$, we have $b = r\text{-lim}_{i \to \infty} f(a_i)$.

**Example 3.5.** Number 2 is a 2-limit and number ½ is a ½-limit of the function $f(x) = x$ at the interval $[0, 1]$.

**Theorem 3.9.** For any point $a \in \mathbf{R}$, $b = (q, r)\text{-lim}_{x \to a} f(x)$ if and only if $b = r\text{-lim}_{x \to [u, v]} f(x))$ with $u = a - q$ and $v = a + q$.

Proof. Necessity. Let $b = (q, r)\text{-lim}_{x \to a} f(x)$ for some point $a$ and $d$ be an arbitrary point from the interval $[a - q, a + q]$. If $d = \lim l$ for a sequence $l = \{a_i \in \text{Dom} f; i \in \omega, a_i \neq a\}$, then results from (Burgin, 2000) imply that $a = q\text{-lim } l$. This, by Definition 3.7, results in $b = r\text{-lim}_{i \to \infty} f(a_i)$. As the point $d$ is arbitrary in the interval $[a - q, a + q]$, we have, by Definition 3.9, that $b = r\text{-lim}_{x \to [u, v]} f(x))$ with $u = a - q$ and $v = a + q$.

Sufficiency. Assume that $b = r\text{-lim}_{x \to [u, v]} f(x))$ with $u = a - q$ and $v = a + q$ and let us consider a sequence $h = \{c_i \in \text{Dom} f; i \in \omega, c_i \neq a\}$ such that $a = q\text{-lim } h$. As a bounded sequence, $h$ is a union of its convergent subsequences $h_t$, $t \in A$. By Lemma 2.2, if $a = q\text{-lim } h$, then $a = q\text{-lim } k$ for any subsequence $k$ of $h$. Let us take a convergent subsequence $h_t = \{c_i \in \text{Dom} f; i \in \omega, c_i \neq a\}$ of the sequence $h = \{c_i \in \text{Dom} f; i \in \omega, c_i \neq a\}$. The condition $a = q\text{-lim } h_t$ implies that $\rho(a, d) < q$ where $d = \lim h_t$ (Burgin, 2000). Consequently, the point $d$ belongs to the interval $[u, v]$ and by our assumption, $b = r\text{-lim}_{i \to \infty} f(d_i)$. This is true for all subsequences of $h$. By Lemma 2.3, $b = r\text{-lim}_{i \to \infty} f(c_i)$. As $h$ is an arbitrary sequence for which the condition $a = q\text{-lim } h$ is valid, we have (by Definition 3.1) that $b = (q, r)\text{-lim}_{x \to a} f(x)$.

Theorem 3.9 is proved.

Properties of fuzzy limits of function at intervals are similar to properties of fuzzy limits of function at points. For instance, Lemma 3.2 implies the following result.

**Lemma 3.7.** For any interval $[u, v]$, if $b = r\text{-lim}_{x \to [u, v]} f(x)$, then $b = q\text{-lim}_{x \to [u, v]} f(x)$ for any $q \geq r$.

**Lemma 3.8.** If $b = r\text{-lim}_{x \to [u, v]} f(x)$ and $[u, v] \subseteq [h, k]$, then $b = r\text{-lim}_{x \to [h, k]} f(x)$.

The concept of an $r$-limit of a function $f$ at an interval gives a new concepts for the classical case.

**Definition 3.10**. A number $b$ is called a limit of a function $f$ at an interval $[u, v]$ (it is denoted by $b = \lim_{x \to [u, v]} f(x)$) if for any point $a$ from the interval $[u, v]$ and any sequence $l = \{a_i \in \text{Dom } f; i \in \omega, a_i \neq a\}$ with $a = \lim l$, we have $b = \lim_{i \to \infty} f(a_i)$.

However, this construction is not very interesting as the following proposition shows.

**Proposition 3.1**. If a limit a function $f(x)$ at an interval $[u, v]$ exists, then $f(x)$ is almost constant on $[u, v]$, that is, for all but finitely many points from $[u, v]$, the values of $f(x)$ coincide.

**Definition 3.11**. A number $b$ is called a weak limit of a function $f$ at an interval $[u, v]$ (it is denoted by $b = \text{wlim}_{x \to [u, v]} f(x)$) if there is a point $a$ from the interval $[u, v]$ such that for any sequence $l = \{a_i \in \text{Dom } f; i \in \omega, a_i \neq a\}$, the condition $a = \lim l$ implies $b = \lim_{i \to \infty} f(a_i)$.

**Proposition 3.2**. The set L( $\text{wlim}_{x \to [u, v]} f(x)$ ) of all weak limits of a function $f$ at an interval $[u, v]$ is equal to $\{ b;\ b = \lim_{x \to a} f(x), a \in [u, v] \}$.

**Remark 3.11.** It is known that a function $f(x)$ is continuous at a point $a \in \mathbf{R}$ if the equality $f(a) = \lim_{x \to a} f(x)$ is true. A similar statement for weak limits at an interval is not true as the following example demonstrates.

**Example 3.6.** Let us consider the following function $f$

$$f(x) = \begin{cases} 1/3 & \text{when } x = 2/3; \\ 2/3 & \text{when } x = 1/3; \\ x & \text{otherwise}. \end{cases}$$

Then $f([0, 1]) = L(\lim_{x \to [0, 1]} f(x))$. However, the function $f(x)$ is not continuous in the interval $[0, 1]$.

**Definition 3.12.** A number $b$ is called a weak $(q, r)$-limit of a function $f$ at a point $a \in \mathbf{R}$ (it is denoted by $b = (q, r)$-wlim$_{x \to a} f(x)$) if there is a sequence $l = \{a_i \in \text{Dom} f ; i \in \omega, a_i \neq a\}$ such that $a = q$-lim $l$ and $b = r$-lim$_{i \to \infty} f(a_i)$.

**Lemma 3.9.** The number $b = r$-wlim$_{x \to a} f(x)$ if and only if $b = (0, r)$-wlim$_{x \to a} f(x)$.

This result demonstrates that the concept of a weak $(q, r)$-limit of a function is a natural extension of the concept of a weak $r$-limit of a function.

**Lemma 3.10.** If $a = (q, r)$-wlim$_{x \to a} f(x)$, then $a = (u, v)$-wlim$_{x \to a} f(x)$ for any $v \geq r$ and $u \geq q$.

Proof is similar to the proof of Lemma 3.6.

Theorem 2.2 from (Burgin, 2000) and Theorem 3.9 imply the following result.

**Theorem 3.10.** For an arbitrary interval $[u, v]$, all $r$-limits at points from the interval $[u, v]$ of a bounded in $[u, v]$ function $f(x)$ belong to some finite interval, the length of which is equal to $v - u + 2r$.

**Remark 3.12.** For unbounded functions, this property can be false.

**Theorem 3.11.** Let $b = (q, r)$-lim$_{x \to a} f(x)$ and $c = (h, k)$-lim$_{x \to a} g(x)$. Then:

a) $b + c = (u, v)$-lim$_{x \to a} (f + g)(x)$ where $v = r + k$ and $u = \min \{q, h\}$;

b) $b - c = (u, v)$-lim$_{x \to a} (f - g)(x)$ where $v = r + k$ and $u = \min \{q, h\}$;

c) $kb = (q, |k| \cdot r)$-lim$_{x \to a} (kf)(x)$ for any $k \in \mathbf{R}$ where $(kf)(x) = k \cdot f(x)$.

Proof. a) Let $b = (q, r)$-lim$_{x \to a} f(x)$ and $c = (h, k)$-lim$_{x \to a} g(x)$. Then by Definition 3.7, for any sequence $l = \{b_i \in \text{Dom} f; i \in \omega, b_i \neq a\}$, the condition $a = q$-lim $l$ implies $b = r$-lim$_{i \to \infty} f(b_i)$ and for any sequence $t = \{c_i \in \text{Dom} g; i \in \omega, c_i \neq a\}$, the condition $a = h$-lim $t$ implies $c = k$-lim$_{i \to \infty} f(c_i)$. Let $u = \min \{q, h\}$. Then for any sequence $d = \{d_i \in \text{Dom} f \cap \text{Dom} g; i \in \omega, d_i \neq a\}$, the condition $a = u$-lim $d$ implies (by Lemma 2.1) the condition $a = q$-lim $d$ and $a = h$-lim $d$. Thus, by Definition 3.7, we have $b = r$-lim$_{i \to \infty} f(d_i)$ and $c = k$-lim$_{i \to \infty} f(d_i)$. By Theorem 2.5 from (Burgin, 2001), we have $b + c = (r + $ $k$)-$\lim_{i \to \infty} f(d_i)$. As $d$ is an arbitrary sequence of elements from the domain Dom $f + g$, $b + c = (u, v)$-$\lim_{x \to a} (f + g)(x)$ by Definition 3.7.

Parts b) and c) are proved in a similar way.

**Proposition 3.1**. If $g(x) \leq f(x) \leq h(x)$ in some neighborhood of the point $a$ that contains the interval $[a - q, a + q]$, $b = (q, r)$-$\lim_{x \to a} g(x)$ and $b = (q, r)$-$\lim_{x \to a} h(x)$ ($b = (q, r)$-wlim$_{x \to a} g(x)$ and $b = (q, r)$-wlim$_{x \to a} h(x)$), then $b = (q, r)$-$\lim_{x \to a} f(x)$ ( $b = (q, r)$-wlim$_{x \to a} f(x)$ ).

**Corollary 3.11**. If $g(x) \leq f(x) \leq h(x)$ in some neighborhood of the point $a$, $b = r$-$\lim_{x \to a} g(x)$ and $b = r$-$\lim_{x \to a} h(x)$ ($b = r$-wlim$_{x \to a} g(x)$ and $b = r$-wlim$_{x \to a} h(x)$), then $b = r$-$\lim_{x \to a} f(x)$ ( $b = r$-wlim$_{x \to a} f(x)$ ).

**Corollary 3.12**. If $g(x) \leq f(x) \leq h(x)$ in some neighborhood of the point $a$, $b = \lim_{x \to a} g(x)$ and $b = \lim_{x \to a} h(x)$ ($b =$ wlim$_{x \to a} g(x)$ and $b =$ wlim$_{x \to a} h(x)$), then $b = \lim_{x \to a} f(x)$ ( $b =$ wlim$_{x \to a} f(x)$ ).

# R e f e r e n c e s